\documentstyle{article}

\title{ZETA FUNCTIONS WITH RESPECT TO GENERAL COINED QUANTUM WALK OF PERIODIC GRAPHS}
\author{Takashi KOMATSU \\
Department of Bioengineering School of Engineering,\\ 
The University of Tokyo \\
Bunkyo, Tokyo, 113-8656, JAPAN \\ 
e-mail: komatsu@coi.t.u-tokyo.ac.jp \\ 
Norio KONNO \\
Department of Applied Mathematics, Faculty of Engineering, \\ 
Yokohama National University \\
Hodogaya, Yokohama 240-8501, JAPAN \\
e-mail: konno-norio-bt@ynu.ac.jp, Tel.: +81-45-339-4205, Fax: +81-45-339-4205 \\ 
Iwao SATO \\ 
Oyama National College of Technology \\
Oyama, Tochigi 323-0806, JAPAN \\ 
e-mail: isato@oyama-ct.ac.jp }
 
\begin{document}
 \maketitle

\clearpage

\vspace{5mm}

{\bf 2000 Mathematical Subject Classification}: 60F05, 05C50, 15A15, 05C25. 

{\bf Key words}: zeta function, quantum walk,  periodic graph, general coin

\vspace{5mm}

The contact author for correspondence: 

Iwao Sato 

Oyama National College of Technology, 
Oyama, Tochigi 323-0806, JAPAN

Tel: +81-285-20-2176

Fax: +81-285-20-2880

E-mail: isato@oyama-ct.ac.jp

\clearpage

\begin{abstract}
We define a zeta function of a graph by using the time evolution matrix of 
a general coined quantum walk on it, and give a determinant expression for 
the zeta function of a finite graph. 
Furthermore, we present a determinant expression for the zeta function an 
(inifinite) periodic graph. 
\end{abstract}

\section{Introduction}

All graphs in this paper are assumed to be simple. 
Let $G$ be a connected graph with vertex set $V(G)$ and edge set $E(G)$, 
and let $R(G)= \{ (u,v),(v,u) \mid uv \in E(G) \} $ be the set of 
oriented edges (or arcs) $(u,v),(v,u)$ directed oppositely for 
each edge $uv$ of $G$. 
For $e=(u,v) \in R(G)$, $u=o(e)$ and $v=t(e)$ are called 
the {\em origin} and the {\em terminal} of $e$, respectively. 
Furthermore, let $e^{-1}=(v,u)$ be the {\em inverse} of $e=(u,v)$. 

A {\em path $P$ of length $n$} in $G$ is a sequence 
$P=(e_1, \cdots ,e_n )$ of $n$ arcs such that $e_i \in R(G)$,
$t( e_i )=o( e_{i+1} )(1 \leq i \leq n-1)$. 
If $e_i =( v_{i-1} , v_i ), \  1 \leq i \leq n$, then we also denote 
$P$ by $( v_0, v_1 , \cdots ,v_n )$. 
Set $ \mid P \mid =n$, $o(P)=o( e_1 )$ and $t(P)=t( e_n )$. 
Also, $P$ is called an {\em $(o(P),t(P))$-path}. 
A $(v, w)$-path is called a {\em $v$-closed path} if $v=w$. 
The {\em inverse} of a closed path $C=( e_1, \cdots ,e_n )$ is the closed 
path $C^{-1} =( e^{-1}_n , \cdots ,e^{-1}_1 )$. 

We say that a path $P=(e_1, \cdots ,e_n )$ has a {\em backtracking} 
or a {\em bump} at $t( e_i )$ 
if $ e^{-1}_{i+1} =e_i $ for some $i(1 \leq i \leq n-1)$. 
A path without backtracking is called {\em proper}. 
Let $B^r$ be the closed path obtained by going $r$ times around a closed 
path $B$. 
Such a closed path is called a {\em multiple} of $B$. 
Multiples of a closed path without bumps may have a bump. 
Such a closed path is said to have a {\em tail}. 
If its length is $n$, then the closed path can be written as 
\[
( e_1 , \cdots , e_k , f_1 , f_2 , \cdots , f_{n-2k} , 
e^{-1}_k , \cdots , e^{-1}_1 ) ,
\]
where $( f_1 , f_2 , \cdots , f_{n-2k} )$ is a closed path. 
A closed path is called {\em reduced} if $C$ has no backtracking nor tail. 
Furthermore, a closed path $C$ is {\em primitive} if it is not a multiple of 
a strictly shorter closed path. 
Let ${\cal C}$ be the set of closed paths. 
Furthermore, let ${\cal C}^{nontail}$ and ${\cal C}^{tail}$ be the set of 
closed paths without tail and closed paths with tail.  
Note that ${\cal C}= {\cal C}^{nontail} \cup {\cal C}^{tail}$ and 
${\cal C}^{nontail} \cap {\cal C}^{tail} = \phi $.

We introduce an equivalence relation between closed paths. 
Two closed paths $C_1 =(e_1, \cdots ,e_m )$ and 
$C_2 =(f_1, \cdots ,f_m )$ are called {\em equivalent} if 
there exists an integer $k$ such that $f_j =e_{j+k} $ for all $j$, 
where the subscripts are read modulo $n$. 
The inverse of $C$ is not equivalent to $C$ if $\mid C \mid \geq 3$. 
Let $[C]$ be the equivalence class which contains a closed path $C$. 
Also, $[C]$ is called a {\em cycle}. 
An equivalence class of primitive closed paths in $G$ is called 
a {\em primitive cycle}. 

Let ${\cal K}$ be the set of cycles of $G$. 
Denote by ${\cal R}$, ${\cal P} \subset {\cal R}$ and 
${\cal PK} \subset {\cal K}$ the set of reduced 
cycles, primitive, reduced cycles and primitive cycles of $G$, respectively. 
Also, primitive, reduced cycles are called {\em prime cycles}. 
Let ${\cal C}_m , {\cal C}^{nontail}_m , {\cal C}^{tail}_m , {\cal K}_m $ 
and ${\cal PK}_m $ be the subset of 
${\cal C}, {\cal C}^{nontail}, {\cal C}^{tail}, {\cal K}$ and ${\cal PK}$ 
consisting of elements with length $m$, respectively. 
Note that each equivalence class of primitive, reduced closed paths 
of a graph $G$ passing through a vertex $v$ of $G$ corresponds 
to a unique conjugacy class of the fundamental group $ \pi {}_1 (G,v)$ 
of $G$ at $v$. 

The {\em Ihara zeta function} of a graph $G$ is 
a function of a complex variable $t$ with $\mid t \mid $ 
sufficiently small, defined by 
\[
{\bf Z} (G, t)= {\bf Z}_G (t)= \prod_{[C] \in {\cal P}} 
(1- t^{ \mid C \mid } )^{-1} ,
\]
where $[C]$ runs over all prime cycles of $G$.

Let $G$ be a connected graph with $n$ vertices $v_1, \cdots ,v_n $. 
The {\em adjacency matrix} ${\bf A}= {\bf A} (G)=(a_{ij} )$ is 
the square matrix such that $a_{ij} =1$ if $v_i$ and $v_j$ are adjacent, 
and $a_{ij} =0$ otherwise.
The {\em degree} of a vertex $v_i$ of $G$ is defined by 
$ \deg v_i = \deg {}_G v_i = \mid \{ v_j \mid v_i v_j \in E(G) \} \mid $. 
If $ \deg {}_G v=k$(constant) for each $v \in V(G)$, then $G$ is called 
{\em $k$-regular}. 

Starting from $p$-adic Selberg zeta functions, Ihara introduced zeta functions 
of graphs [13]. 
At the beginning, Serre [20] pointed out that the Ihara zeta function is 
the zeta function of a regular graph. 
Ihara [13] showed that the reciprocal of the Ihara zeta function of a regular graph 
is an explicit polynomial. 

The Ihara zeta function of a regular graph has the following three properties: 
rationality; the functional equations; 
the analogue of the Riemann hypothesis(see [24]).  
The analogue of the Riemann hypothesis for the Ihara zeta function of a graph is 
given as follows: 
Let $G$ be any connected $(q+1)$-regular graph$(q>1)$ and 
$s= \sigma +it \  ( \sigma ,t \in {\bf R})$ a complex number. 
If ${\bf Z}_G ( q^{-s} )^{-1} =0$ and ${\rm Re} \  s \in (0,1)$, then 
${\rm Re} \  s = \frac{1}{2} $. 
A connected $(q+1)$-regular graph $G$ is called a 
{\em Ramanujan graph} if for all eigenvalues $\lambda $ of the 
adjacency matrix ${\bf A} (G)$ of $G$ such that 
$ \lambda \neq \pm (q+1)$, we have $\mid \lambda \mid \leq 2 \sqrt{q} $. 
This definition was introduced by Lubotzky, Phillips and Sarnak [16]. 
For a connected $(q+1)$-regular graph $G$, ${\bf Z}_G ( q^{-s} )$ satisfies 
the Riemann hypothesis if and only if $G$ is a Ramanujan graph. 

A zeta function of a regular graph $G$ associated to a unitary 
representation of the fundamental group of $G$ was developed by 
Sunada [22,23]. 
Hashimoto [11] treated multivariable zeta functions of bipartite graphs. 
Bass [1] generalized Ihara's result on the Ihara zeta function of 
a regular graph to an irregular graph, and showed that its reciprocal is 
a polynomial.

\newtheorem{theorem}{Theorem}
\begin{theorem}[Bass] 
Let $G$ be a connected graph. 
Then the reciprocal of the Ihara zeta function of $G$ is given by 
\[
{\bf Z} (G,t )^{-1} =(1- t^2 )^{r-1} 
\det ( {\bf I} -t {\bf A} (G)+ t^2 ( {\bf D} - {\bf I} )) , 
\]
where $r$ is the Betti number of $G$, 
and ${\bf D} =( d_{ij} )$ is the diagonal matrix 
with $d_{ii} = \deg v_i$ and $d_{ij} =0, i \neq j , 
(V(G)= \{ v_1 , \cdots , v_n \} )$. 
\end{theorem}

Stark and Terras [21] gave an elementary proof of Theorem 1, and 
discussed three different zeta functions of any graph. 
Various proofs of Bass' theorem are known. 
Kotani and Sunada [15] proved Bass' theorem by using the property of the 
Perron operator. 
Foata and Zeilberger [6] presented a new proof of Bass' theorem by 
using the algebra of Lyndon words.

The Ihara zeta function of a finite graph was extended to an infinite graph 
in [1,3,7,8,9,10], and its determinant expressions were presented. 
Bass [2] defined the zeta function for a pair of a tree $X$ and a countable 
group $\Gamma $ which acts discretely on $X$ with quotient being a graph of 
finite groups. 
Clair and Mokhtari-Sharghi [3] extended Ihara zeta functions to infinite 
graphs on which a group $\Gamma $ acts isomorphically and with finite 
quotient. 
In [7], Grigorchuk and \.{Z}uk defined zeta functions of infinite discrete 
groups, and of some class of infinite periodic graphs. 

Guido, Isola and Lapidus [8] defined the Ihara zeta function of a 
periodic simple graph. 
Let $G=(V(G),E(G))$ be a countable simple graph, 
and let $\Gamma $ be a countable discrete subgroup of 
automorphisms of $G$, which acts freely on $G$, and with finite quotient 
$G/ \Gamma $. 
Then the Ihara zeta function of a periodic simple graph is defined 
as follows: 
\[
{\bf Z}_{G, \Gamma} (t)= \prod_{[C]_{\Gamma} \in [{\cal P}]_{\Gamma }} 
(1- t^{\mid C \mid } )^{- 1/ \mid \Gamma {}_{[C]} \mid } , 
\]
where $ \Gamma {}_{[C]} $ is the stabilizer of $[C]$ in $\Gamma $, and 
$[C]_{\Gamma} $ runs over all $\Gamma $-equivalence classes of 
prime cycles in $G$. 

Guido, Isola and Lapidus [8] presented a determinant expression for 
the Ihara zeta function of a periodic simple graph by using Stark and Terras' 
method [21].

\begin{theorem}[Guido, Isola and Lapidus] 
For a periodic simple graph $G$, 
\[
{\bf Z}_{G, \Gamma } (t)=(1- t^2 )^{-(m-n)} 
\det {}_{\Gamma } ( {\bf I} -t {\bf A} (G) +( {\bf D} - {\bf I} ) 
t^2 )^{-1} , 
\]
where $\det {}_{\Gamma } $ is a determinant for bounded operators 
belonging to a von Neumann algebra with a finite trace. 
\end{theorem}

Guido, Isola and Lapidus [9] presented a determinant expression for 
the Ihara zeta function of a periodic graph $G$ and a countable discrete subgroup 
$\Gamma $ of aoutomorphisms of $G$ which acts discretely without inversions, and 
with bounded covolume.

\begin{theorem}[Guido, Isola and Lapidus] 
For a periodic graph $G$, 
\[
{\bf Z}_{G, \Gamma } (t)^{-1} =(1- t^2 )^{\chi {}^{(2)} (G)} 
\det {}_{\Gamma } ( \Delta (t)) , 
\]
where ${\chi {}^{(2)} (G)}$ is the $L^2$-Euler characteristic of $(G, \Gamma )$
(see [3]), and $\Delta (t)={\bf I} -t{\bf A}+ t^2 ( {\bf D} -{\bf I} )$. 
\end{theorem}

\section{The Grover walk on a graph}

A discrete-time quantum walk is a quantum analog of the classical random walk on a graph whose state vector is governed by 
a matrix called the transition matrix.  
Let $G$ be a connected graph with $n$ vertices and $m$ edges, 
$V(G)= \{ v_1 , \ldots , v_n \} $ and $R(G)= \{ e_1 , \ldots , e_m , 
e^{-1}_1 , \ldots , e^{-1}_m \} $. 
Set $d_j = d_{u_j} = \deg v_j $ for $i=1, \ldots ,n$. 
The {\em transition matrix} ${\bf U} ={\bf U} (G)=( U_{ef} )_{e,f \in R(G)} $ 
of $G$ is defined by 
\[
U_{ef} =\left\{
\begin{array}{ll}
2/d_{t(f)} (=2/d_{o(e)} ) & \mbox{if $t(f)=o(e)$ and $f \neq e^{-1} $, } \\
2/d_{t(f)} -1 & \mbox{if $f= e^{-1} $, } \\
0 & \mbox{otherwise.}
\end{array}
\right. 
\]
The matrix ${\bf U} $ is called the {\em Grover matrix} of $G$. 

We introduce the {\em positive support} $\>{\bf F}^+ =( F^+_{ij} )$ of 
a real matrix ${\bf F} =( F_{ij} )$ as follows: 
\[
F^+_{ij} =\left\{
\begin{array}{ll}
1 & \mbox{if $F_{ij} >0$, } \\
0 & \mbox{otherwise.}
\end{array}
\right.
\]

Let $G$ be a connected graph. 
If the degree of each vertex of $G$ is not less than 2, i.e., $\delta (G) \geq 2$, 
then $G$ is called an {\em md2 graph}. 

The transition matrix of a discrete-time quantum walk in a graph 
is closely related to the Ihara zeta function of a graph. 
We stare a relationship between the discrete-time quantum walk and the Ihara zeta function of a graph by Ren et al. [19].

Konno and Sato [14] obtained the following formula of the characteristic polynomial of ${\bf U}$ 
by using the determinant expression for the second weighted zeta function of a graph. 

Let $G$ be a connected graph with $n$ vertices and $m$ edges. 
Then the $n \times n$ matrix ${\bf T} (G)=( T_{uv} )_{u,v \in V(G)}$ is given as follows: 
\[
T_{uv} =\left\{
\begin{array}{ll}
1/( \deg {}_G u)  & \mbox{if $(u,v) \in R(G)$, } \\
0 & \mbox{otherwise.}
\end{array}
\right.
\] 
Note that the matrix ${\bf T} (G)$ is the transition matrix of the simple random walk on $G$.

\begin{theorem}[Konno and Sato]
Let $G$ be a connected graph with $n$ vertices $v_1 , \ldots , v_n $ and $m$ edges. Then, for the transition matrix ${\bf U}$ of $G$, we have
\[
\begin{array}{rcl}
\det ( \lambda {\bf I}_{2m} - {\bf U} ) & = & ( \lambda {}^2 -1)^{m-n} \det (( \lambda {}^2 +1) {\bf I}_n -2 \lambda {\bf T} (G)) \\
\  &   &                \\ 
\  & = & \frac{( \lambda {}^2 -1)^{m-n} \det (( \lambda {}^2 +1) {\bf D} -2 \lambda {\bf A} (G))}{d_{v_1} \cdots d_{v_n }} .   
\end{array}
\]
\end{theorem}

From this Theorem, the spectra of the Grover matrix on a graph is obtained by means of those of 
${\bf T} (G)$ (see [4]). 
Let $Spec ({\bf F})$ be the spectra of a square matrix ${\bf F}$ .

\newtheorem{corollary}{Corollary} 
\begin{corollary}[Emms, Hancock, Severini and Wilson] 
Let $G$ be a connected graph with $n$ vertices and $m$ edges. 
The transition matrix ${\bf U}$ has $2n$ eigenvalues of the form 
\[
\lambda = \lambda {}_T \pm i \sqrt{1- \lambda {}^2_T } , 
\]
where $\lambda {}_T $ is an eigenvalue of the matrix ${\bf T} (G)$. 
The remaining $2(m-n)$ eigenvalues of ${\bf U}$ are $\pm 1$ with equal multiplicities. 
\end{corollary}

In this paper, we define a zeta function of a periodic graph by using the time evolution matrix 
of a general coined quantum walk on it, and present its determinant expression. 
The proof is an analogue of Bass' method [1].  

In Section 1, we state a review for the Ihara zeta function of a finite graph and 
infinite graphs, i.e., a periodic simple graph, a periodic graph. 
In Section 2, we state about the Grover walk on a graph as a discrete-time quantum walk on a graph. 
In Section 3, we define a zeta function of a finite graph $G$ by using the time evolution matrix 
of a general coined quantum walk on $G$, and present its determinant expression. 
Furthermore, we give the explicit formula for the characteristic polynomial of the time evolution matrix 
of a general coined quantum walk on $G$, and so present its spectrum. 
In Section 4, we state the definition of a periodic graph. 
In Section 5, we review a determinant for bounded operators acting on an 
infinite dimensional Hilbert space and belonging to a von Neumann algebra 
with a finite trace. 
In Section 6, we present a determinant expression for the above zeta function of a periodic graph.

\section{Spectra for the unitary matrix of a general coined quantum walk on a graph} 

We consider a generalization of a coined quantum walk on a graph. 
We replace the coin operator ${\bf C} $ of a coined quantum walk with unitary matrix with two spectra which are distinct 
from $\pm 1$. 

For given connected graph $G$ with $n$ vertices and $m$ edges, let ${\bf d} : \ell {}^2 (V(G)) \longrightarrow \ell {}^2 (R(G))$ 
such that 
\[
{\bf d} {\bf d}^* ={\bf I}_{n} . 
\] 
Furthermore, let 
\[
{\bf C} = a {\bf d}^* {\bf d}+b( {\bf I}_{2m} - {\bf d}^* {\bf d} ) 
\] 
and ${\bf U} = {\bf S} {\bf C} $(see [12]). 
A discrete-time quantum walk on $G$ with ${\bf U}$ as a time evolution matrix is called a {\em general coined quantum walk} 
on $G$. 
Then we define a zeta function of $G$ by using ${\bf U}$ as follows: 
\[
\zeta {} (G, u)= \det ( {\bf I}_{2m} -u {\bf U} )^{-1} 
= \det ( {\bf I}_{2m} -u {\bf S} (a {\bf d}^* {\bf d} +b( {\bf I}_{2m} - {\bf d}^* {\bf d} )) )^{-1} . 
\]

Now, we have the following result.

\begin{theorem}
Let $G$ be a connected graph $n$ vertices and $m$ edges, ${\bf U} = {\bf S} {\bf C} $ the time evolution matrix of a coined quantum walk on $G$. 
Suppose that $\sigma ( {\bf C} )= \{ a,b \} $. 
Set $q= \dim \ ker(a-C)$.
Then, for the unitary matrix ${\bf U} = {\bf S} {\bf C} $, we have    
\[
\zeta {} (G, u)=(1- b^2 u )^{m-n} \det ((1-abu^2 ) {\bf I}_n -cu {\bf d} {\bf S} {\bf d}^* ) , c=a-b .   
\]
\end{theorem}

{\em Proof }.  At first, we have  
\[
\begin{array}{rcl}
\  &   & \zeta {} (G, u)= \det ( {\bf I}_{2m} -u {\bf U} )= \det ( {\bf I}_{2m} -u {\bf S} {\bf C} ) \\ 
\  &   &                \\ 
\  & = & 
\det ({\bf I}_{2m} -u {\bf S} (a {\bf d}^* {\bf d} +b( {\bf I}_{2m} - {\bf d}^* {\bf d} ))) \\ 
\  &   &                \\ 
\  & = & 
\det ({\bf I}_{2m} -u {\bf S} ((a-b) {\bf d}^* {\bf d} +b {\bf I}_{2m} ))) \\ 
\  &   &                \\ 
\  & = & \det ( {\bf I}_{2m} -bu {\bf S} -cu {\bf S} {\bf d}^* {\bf d} ) \\ 
\  &   &                \\
\  & = & \det ( {\bf I}_{2m} - cu {\bf S} {\bf d}^* {\bf d} ( {\bf I}_{2m} -bu {\bf S} )^{-1} ) \det ( {\bf I}_{2m} -bu {\bf S} ) . 
\end{array}
\]

But, if ${\bf A}$ and ${\bf B}$ are a $m \times n $ and $n \times m$ 
matrices, respectively, then we have 
\[
\det ( {\bf I}_{m} - {\bf A} {\bf B} )= 
\det ( {\bf I}_n - {\bf B} {\bf A} ) . 
\]
Thus, we have 
\[
\det ( {\bf I}_{2m} -u {\bf U} )= \det ( {\bf I}_{2m} -u {\bf S} {\bf C} )= 
\det ( {\bf I}_{n} - cu {\bf d} ( {\bf I}_{2m} -bu {\bf S} )^{-1} {\bf S} {\bf d}^* ) \det ( {\bf I}_{2m} -bu {\bf S} ) . 
\]
But, we have 
\[
\det ( {\bf I}_{2m} -bu {\bf S} )=(1-b^2 u^2 )^m . 
\]
Furthermore, we have 
\[
( {\bf I}_{2m} -bu {\bf S} )^{-1} = \frac{1}{1-b^2 u^2 } ({\bf I}_{2m} +u {\bf S} ) . 
\]

Therefore, it follows that 
\[
\begin{array}{rcl}
\  &   & \det ( {\bf I}_{2m} -u {\bf U} ) \\ 
\  &   &                \\ 
\  & = & (1-b^2 u^2 )^{m} 
\det ( {\bf I}_{2m} - \frac{cu}{1-b^2 u^2 } {\bf S} {\bf d}^* {\bf d} ( {\bf I}_{2m} +bu {\bf S} )) \\ 
\  &   &                \\ 
\  & = & (1-b^2 u^2 )^{m} 
\det ( {\bf I}_{n} - \frac{cu}{1-b^2 u^2} {\bf d} ( {\bf I}_{2m} +bu {\bf S} ) {\bf S} {\bf d}^* ) \\ 
\  &   &                \\ 
\  & = & (1-b^2 u^2 )^{m-n} 
\det ((1- b^2 u^2 ) {\bf I}_{n} -cu {\bf d} {\bf S} {\bf d}^* -bc u^2 {\bf d} {\bf S}^2 {\bf d}^* ) \\ 
\  &   &                \\ 
\  & = & (1-b^2 u^2 )^{m-n} 
\det ((1- b^2 u^2 ) {\bf I}_{n} -cu {\bf d} {\bf S} {\bf d}^* -bc u^2 {\bf I}_n ) \\ 
\  &   &                \\ 
\  & = & (1-b^2 u^2 )^{m-n} 
\det ((1-ab u^2 ) {\bf I}_{n} -cu {\bf d} {\bf S} {\bf d}^* ) . 
\end{array}
\]
$\Box$

\begin{corollary}
Let $G$ be a connected with $n$ vertices and $m$ edges. 
Then, for the unitary matrix ${\bf U} = {\bf S} {\bf C} $, we have    
\[
\det ( \lambda {\bf I}_{2m} -u {\bf U} )=(\lambda {}^2 -b^2 )^{m-q} \det (( \lambda {}^2 -ab) {\bf I}_q -c \lambda {\bf d} {\bf S} {\bf d}^* ) ,  
\]
where $q= \dim \ ker(1- {\bf C} )$. 
\end{corollary}

{\em Proof }.  Let $u=1/ \lambda $. 
Then, by Theorem 1.1, we have 
\[
\det ({\bf I}_{2m} -1/  \lambda {\bf U} )=(1-b^2 / \lambda {}^2 )^{m-q} \det ((1-ab/ \lambda {}^2 ) {\bf I}_q - c/\lambda {\bf d} {\bf S} {\bf d}^* ) ,  
\]
and so, 
\[
\det ( \lambda  {\bf I}_{2m} - {\bf U} )=( \lambda {}^2 -b^2 )^{m-q} \det (( \lambda {}^2 -ab) {\bf I}_q -c \lambda {\bf d} {\bf S} {\bf d}^* ) . 
\]
$\Box$

By Corollary 3.2, the following result holds.

\begin{corollary}
Let $G$ be a connected with $n$ vertices and $m$ edges. 
Then, the spectra of the unitary matrix ${\bf U} = {\bf S} {\bf C} $ are given as follows: 
\begin{enumerate} 
\item $2q$ eigenvalues: 
\[
\lambda = \frac{c \mu \pm \sqrt{ c^2 \mu {}^2 +4ab}}{2} , \ \mu \in Spec ( {\bf d} {\bf S} {\bf d}^* ) ; 
\]
\item The rest eigenvalues are $\pm b$ with the same multiplicity $m-q$. 
\end{enumerate} 
\end{corollary}

{\em Proof }. By Corollary 3.2, we have 
\[
\begin{array}{rcl}
\  &   & \det ( \lambda  {\bf I}_{2m} - {\bf U} ) \\
\  &   &                \\ 
\  & = & ( \lambda {}^2 -b^2 )^{m-q} 
\prod_{ \mu \in Spec ( {\bf d} {\bf S} {\bf d}^* )} ( \lambda {}^2 -c \mu \lambda -ab) . 
\end{array}
\]

Solving $ \lambda {}^2 -2 \mu \lambda +1=0$, we obtain 
\[
\lambda = \frac{c \mu \pm \sqrt{ c^2 \mu {}^2 +4ab}}{2} . 
\] 
The result follows. 
$\Box$

\section{Periodic graphs}

Let $G=(V(G),E(G))$ be a simple graph. 
Assume that $G$ is countable ($V(G)$ and $E(G)$ are countable), 
and with bounded degree, i.e., $d= \sup_{v \in V(G)} \deg v < \infty$. 
Let $\Gamma $ be a countable discrete subgroup of automorphisms of $G$, 
which acts 
\begin{enumerate} 
\item  without inversions: $\gamma (e) \neq e^{-1} $ for any 
$\gamma \in \Gamma, e \in R(G)$, 
\item  discretely: $\Gamma {}_v = \{ \gamma \in \Gamma \mid \gamma v=v \}$ 
is finite for any $v \in V(G)$, 
\item  with bounded covolume: ${\rm vol} (G/ \Gamma ):= 
\sum_{v \in {\cal F}_0 } \frac{1}{ \mid \Gamma {}_v \mid } < \infty $, where 
${\cal F}_0 \subset V(G)$ contains exactly one representative for 
each equivalence class in $V(G/ \Gamma )$. 
\end{enumerate} 
Then $G$ is called a {\em periodic graph} with a countable discrete subgroup 
$\Gamma $ of $Aut \  G$. 
Note that the third condition is equivalent to the following condition: 
\[
{\rm vol} (R(G)/ \Gamma ):= 
\sum_{e \in {\cal F}_1 } \frac{1}{ \mid \Gamma {}_e \mid } < \infty , 
\]
where ${\cal F}_1 \subset R(G)$ contains exactly one representative for 
each equivalence class in $R(G/ \Gamma )$. 

Let $\ell {}^2 (V(G))$ be the Hilbert space of functions 
$f: V(G) \longrightarrow {\bf C} $ such that 
$\mid \mid f \mid \mid := \sum_{v \in V(G)} \mid f(v) \mid {}^2 
< \infty $. 
We define the left regular representation $\lambda {}_0 $ of $\Gamma $ on 
$\ell {}^2 (V(G))$ as follows: 
\[
(\lambda {}_0 ( \gamma )f)(x)=f( \gamma {}^{-1} x), \  \gamma \in \Gamma, 
\  f \in \ell {}^2 (V(G)), \  x \in V(G) . 
\]

We state the definition of a von Neumann algebra. 
Let $H$ be a separable complex Hilbert space, and let ${\cal B} (H)$ denote 
the ${\bf C}^*$-algebra of bounded linear operators on $H$. 
For a subset $M \subset {\cal B} (H)$, the {\em commutant} of $M$ is 
$M^{\prime} = \{ T \in {\cal B} (H) \mid ST=TS, \forall S \in M \} $. 
Then a {\em von Neumann algebra} is a subalgebra ${\cal A} \leq {\cal B} (H)$ 
such that ${\cal A}^{\prime \prime} ={\cal A} $. 
It is known that a determinant is defined for a suitable class of operators 
in a von Neumann algebra with a finite trace (see [5,8]). 
 
For the Hilbert space $\ell {}^2 (V(G))$, we consider a von Neumann algebra. 
Let ${\cal B} (\ell {}^2 (V(G)))$ be the ${\bf C}^*$-algebra of bounded 
linear operators on $\ell {}^2 (V(G))$. 
A bounded linear operator $A$ of ${\cal B} (\ell {}^2 (V(G)))$ acts on 
$\ell {}^2 (V(G))$ by 
\[
A(f)(v)= \sum_{w \in V(G)} A(v,w)f(w), \   v \in V(G), \  f \in \ell {}^2 (V(G)) . 
\]
Then the von Neumann algebra ${\cal N}_0 (G, \Gamma )$ of bounded operators on 
$\ell {}^2 (V(G))$ commuting with the action of $\Gamma $ is defined as follows: 
\[
{\cal N}_0 (G, \Gamma )= \{ \lambda {}_0 (\gamma ) \mid \gamma \in \Gamma \} {}^{\prime} 
= \{ T \in {\cal B} (\ell {}^2 (V(G))) \mid \lambda {}_0 (\gamma )T=T \lambda {}_0 (\gamma ), 
\forall \gamma \in \Gamma \} . 
\]
The von Neumann algebra ${\cal N}_0 (G, \Gamma )$ inherits a trace 
by 
\[
{\rm Tr}_{\Gamma} (A)=\sum_{x \in {\cal F}_0 } \frac{1}{ \mid \Gamma {}_x \mid } A(x,x), 
\  A \in {\cal N}_0 (G, \Gamma ) . 
\]
Let the adjacency matrix ${\bf A} = {\bf A} (G)$ of $G$ be defined by 
\[
( {\bf A} f)(v)= \sum_{(v,w) \in R(G)} f(w), \  f \in \ell {}^2 (V(G)) . 
\]
By [16,17], we have 
\[
\mid \mid {\bf A} \mid \mid \leq d= \sup_{v \in V(G)} \deg {}_G v < \infty , 
\]
and so ${\bf A} \in {\cal N}_0 (G, \Gamma )$.

Similarly to $\ell {}^2 (V(G))$, we consider the Hilbert space 
$\ell {}^2 (R(G))$ of functions $f: R(G) \longrightarrow {\bf C} $ such that 
$\mid \mid \omega \mid \mid := \sum_{e \in R(G)} \mid \omega (e) \mid {}^2 
< \infty $. 
We define the left regular representation $\lambda {}_1 $ of $\Gamma $ on 
$\ell {}^2 (R(G))$ as follows: 
\[
(\lambda {}_1 ( \gamma ) \omega )(e)= \omega ( \gamma {}^{-1} e), 
\  \gamma \in \Gamma, \  \omega \in \ell {}^2 (R(G)), \  e \in R(G) . 
\]
Then the von Neumann algebra ${\cal N}_1 (G, \Gamma )= \{ \lambda {}_1 
(\gamma ) \mid \gamma \in \Gamma \} {}^{\prime} $ of bounded operators on 
$\ell {}^2 (R(G))$ commuting with the action of $\Gamma $, inherits a trace 
by 
\[
{\rm Tr}_{\Gamma} (A)=\sum_{e \in {\cal F}_1 } \frac{1}{ \mid \Gamma {}_e \mid } A(e,e), 
\  A \in {\cal N}_1 (G, \Gamma ) . 
\]

\section{An analytic determinant for von Neumann algebras 
with a finite trace}

In an excellent paper [5], Fuglede and Kadison defined a positive-valued 
determinant for a von Neumann algebra with trivial center and finite trace $\tau $. 
For an invertible operator $A$ with polar decomposition $A=UH$, the 
Fuglede-Kadison determinant of $A$ is defined by 
\[
Det(A)= \exp \circ \tau \circ \log H, 
\]
where $\log H$ may be defined via functional calculus. 

Guido, Isola and Lapidus [8] extended the Fuglede-Kadison determinant 
to a determinant which is an analytic function. 
Let $( {\cal A} , \tau )$ be a von Neumann algebra with a finite trace 
$\tau $. 
Then, for $A \in {\cal A} $, let 
\[
\det {}_{\tau} (A)= \exp \circ \tau \circ \log A, 
\]
where 
\[
\log (A):= \frac{1}{2 \pi i} \int_{\Lambda } \log \lambda 
( \lambda -A )^{-1} d \lambda , 
\]
and $\Lambda $ is the boundary of a connected, simply connected region 
$\Omega $ containing the spectrum $\sigma (A)$ of $A$. 
Then the following lemma holds (see [8, Lemma 5.1]).

\newtheorem{lemma}{Lemma}
\begin{lemma}[Guido, Isola and Lapidus] 
Let ${\cal A}, \Omega, \Gamma $ be as above, and $\phi, \psi$ two branches 
of the logarithm such that both domains contain $\Omega $. 
Then 
\[
\exp \circ \tau \circ \phi (A)= \exp \circ \tau \circ \psi (A) . 
\]
\end{lemma}

Next, we consider a determinant on some subset of ${\cal A}$. 
Let $( {\cal A} , \tau )$ be a von Neumann algebra with a finite trace, and 
${\cal A} {}_0 = \{ A \in {\cal A} \mid 0 \notin {\rm conv} \  \sigma (A) \} $, 
where ${\rm conv} \  \sigma (A)$ is the convex hull of $ \sigma (A)$. 
For any $A \in {\cal A} {}_0 $, we set 
\[
\det {}_{\tau} (A)= \exp \circ \tau \circ 
( \frac{1}{2 \pi i} \int_{\Lambda } \log \lambda 
( \lambda -A )^{-1} d \lambda ) , 
\]
where $\Lambda $ is the boundary of a connected, simply connected region 
$\Omega $ containing the spectrum ${\rm conv} \  \sigma (A)$, and 
$\log $ is a branch of the logarithm whose domain contains $\Omega $. 
Then the above determinant is well-defined and analytic on 
${\cal A} {}_0 $ (see [8, Corollary 5.3]). 
Furthermore, Guido, Isola and Lapidus [8,9] showed that $\det {}_{\tau } $ 
has the following properties.

\newtheorem{proposition}{Proposition}
\begin{proposition}[Guido, Isola and Lapidus] 
Let $( {\cal A} , \tau )$ be a von Neumann algebra with a finite trace, 
$A \in {\cal A} {}_0 $. 
Then 
\begin{enumerate} 
\item  $\det {}_{\tau} (zA)=z^{\tau (I)} \det {}_{\tau} (A)$ for any 
$z \in {\bf C} \setminus \{ 0 \} $. 
\item  If $A$ is normal, and $A=UH$ is its polar decomposition, 
\[
\det {}_{\tau} (A)= \det {}_{\tau} (U) \det {}_{\tau} (H) . 
\]
\item  If $A$ is positive, $\det {}_{\tau} (A)= Det (A)$, 
where $Det (A) $ is the Fuglede-Kadison determinant of $A$. 
\end{enumerate} 
\end{proposition}

\begin{proposition}[Guido, Isola and Lapidus] 
Let $( {\cal A} , \tau )$ be a von Neumann algebra with a finite trace. 
Then 
\begin{enumerate} 
\item  For $A, B \in {\cal A}$ and sufficiently small $u \in {\bf C} $, 
\[
\det {}_{\tau } ((I+uA)(I+uB))= \det {}_{\tau } (I+uA) 
\det {}_{\tau } (I+uB) . 
\]
\item  If $A \in {\cal A}$ has a bounded inverse, and $T \in {\cal A}_0 $, 
then 
\[
\det {}_{\tau } ( AT A^{-1} )= \det {}_{\tau } (T) . 
\]
\item  If 
\[
T= 
\left[
\begin{array}{cc}
T_{11} & T_{12} \\
0 & T_{22} 
\end{array}
\right]
\in {\rm Mat}_2 ( {\cal A} )
,
\]
with $T_{ii} \in {\cal A} $ such that 
$\sigma ( T_{ii} ) \subset B_1 (1) 
:= \{ z \in {\bf C} \mid \  \mid z-1 \mid <1 \}$ for $i=1,2$, then 
\[
\det {}_{\tau } (T)= \det {}_{\tau } ( T_{11} ) 
\det {}_{\tau } ( T_{22} ) . 
\]
\end{enumerate} 
\end{proposition}

\begin{corollary}[Guido, Isola and Lapidus] 
Let $ \Gamma $ be a discrete group , $\pi {}_1 , \pi {}_2 $ unitary 
representations of $\Gamma $, $\tau {}_1 $, $\tau {}_2 $ finite traces on 
$\pi {}_1 ( \Gamma )^{\prime } $ and $\pi {}_2 ( \Gamma )^{\prime } $, respectively. 
Let $\pi =\pi {}_1 \oplus \pi {}_2 $, $\pi =\tau {}_1 + \tau {}_2 $ and 
$T= 
\left[
\begin{array}{cc}
T_{11} & T_{12} \\
0 & T_{22} 
\end{array}
\right]
\in \pi ( \Gamma )^{\prime} $, with 
$\sigma ( T_{ii} ) \subset B_1 (1) 
:= \{ z \in {\bf C} \mid \  \mid z-1 \mid <1 \}$ for $i=1,2$, then 
\[
\det {}_{\tau } (T)= \det {}_{\tau {}_1 } ( T_{11} ) 
\det {}_{\tau {}_2 } ( T_{22} ) . 
\]
\end{corollary}

\section{A zeta function with respect to a general coined quantum walk of an infinite periodic graph}

We define a zeta function with respect to a general coined quantum walk of an infinite periodic graph. 

Let $G$ be a periodic graph with a countable discrete subgroup $\Gamma $ of $Aut \  G$. 
Moreover, let 
\[
{\bf I}_V =Id_{\ell {}^2 (V(G))} , {\bf I}_R =Id_{\ell {}^2 (R(G))} . 
\] 
Then, let ${\bf d} : \ell {}^2 (V(G)) \longrightarrow \ell {}^2 (R(G))$ such that 
\[
{\bf d} {\bf d}^* ={\bf I}_V . 
\] 
Furthermore, let 
\[
{\bf C} = a {\bf d}^* {\bf d} +b( {\bf I}_R - {\bf d}^* {\bf d} ) 
\] 
and ${\bf U} = {\bf S} {\bf C} $, 
where {\bf S} is the operator on $\ell {}^2 (R(G))$ such that 
\[
({\bf S} \omega )(e)= \omega ( e^{-1} ) , \ \omega \in \ell {}^2 (R(G)), \ e \in R(G) .  
\]

Then a zeta function with respect to a general coined quantum walk of $G$ is defined 
as follows: 
\[
\zeta {} (G, \Gamma, u)= \det {}_{\Gamma } ( {\bf I}_R -u {\bf U} )^{-1} 
= \det {}_{\Gamma } ( {\bf I}_R -u {\bf S} (a {\bf d}^* {\bf d} +b( {\bf I}_R - {\bf d}^* {\bf d} )) )^{-1} . 
\]

Then we have the following result.

\begin{theorem}
Let $G$ be a periodic graph with a countable discrete subgroup $\Gamma $ 
of $Aut \  G$. 
Then 
\[
\zeta (G, \Gamma , u)=(1- b^2 u^2 )^{ {\rm Tr}_{\Gamma } ({\bf I}_V)- \frac{1}{2} {\rm Tr}_{\Gamma } ({\bf I}_R)} 
\det {}_{\Gamma } (  (1-ab u^2 ) {\bf I}_V -cu {\bf d} {\bf S} {\bf d}^* ) ,   
\]
where ${\rm Tr}_{\Gamma } ({\bf I}_R)= \sum_{e \in {\cal F}_1 } \frac{1}{ \mid \Gamma {}_e \mid } $ and 
$ {\rm Tr}_{\Gamma } ({\bf I}_V)= \sum_{v \in {\cal F}_0 } \frac{1}{ \mid \Gamma {}_v \mid } $(see [2]). 
\end{theorem}

{\bf Proof}.  The argument is an analogue of the method of Bass [1]. 

Let $G$ be a periodic graph with a countable discrete subgroup $\Gamma $ 
of $Aut \  G$. 

Now we consider the direct sum of the unitary representations 
$\lambda {}_0 $ and $\lambda {}_1 $: 
$\lambda (\gamma ):= \lambda {}_0 (\gamma ) \oplus \lambda {}_1 (\gamma )$  
$ \in {\cal B} (\ell {}^2 (V(G)) \oplus \ell {}^2 (R(G)))$. 
Then the von Neumann algebra $\lambda (\Gamma )^{\prime} 
:=\{ S \in {\cal B} (\ell {}^2 (V(G)) \oplus \ell {}^2 (R(G))) \mid 
S \lambda ( \gamma )= \lambda (\gamma )S, \gamma \in \Gamma \}$ 
consists of operators 
\[
S= 
\left[
\begin{array}{cc}
S_{00} & S_{01} \\
S_{10} & S_{11} 
\end{array}
\right]
,
\]
where $S_{ij} \lambda {}_j (\gamma )=\lambda {}_i (\gamma )S_{ij} , 
\gamma \in \Gamma , i,j=0,1$, so that $S_{ii} \in \Lambda {}_i 
\equiv  {\cal N}_i (G, \Gamma ), i=0,1$. 
Thus, $\lambda ( \Gamma )^{\prime} $ inherits a trace given by 
\[
{\rm Tr }_{\Gamma } 
\left[
\begin{array}{cc}
S_{00} & S_{01} \\
S_{10} & S_{11} 
\end{array}
\right]
:= {\rm Tr}_{\Gamma } (S_{00} )+ {\rm Tr}_{\Gamma } (S_{11} ) . 
\]

We introduce two operators as follows: 
\[
{\bf L} =
\left[
\begin{array}{cc}
(1- b^2 u^2 ) {\bf I}_V & -c {\bf d} -bc u {\bf d} {\bf S} \\
0 & {\bf I}_R 
\end{array}
\right]
,
{\bf M} = 
\left[
\begin{array}{cc}
{\bf I}_V & c {\bf d} +bcu {\bf d} {\bf S} \\
u {\bf S} {\bf d}^* & (1-b^2 u^2 ) {\bf I}_R 
\end{array}
\right]
,  
\]  
where $c=a-b$. 
By 3 and 4 of Lemma 3, we have 
\[
\begin{array}{rcl}
{\bf LM} & = & 
\left[
\begin{array}{cc}
(1- b^2 u^2 ) {\bf I}_V -cu {\bf d} {\bf S} d^* -bc u^2 {\bf d} {\bf S} ^2 {\bf d}^* & 0 \\
u {\bf S} {\bf d}^* & (1- b^2 u^2 ) {\bf I}_R 
\end{array}
\right]
\\ 
\  &   &                \\ 
\  & = & 
\left[
\begin{array}{cc}
(1-ab u^2 ) {\bf I}_V -cu {\bf d} {\bf S} {\bf d}^* & 0 \\
u {\bf S} {\bf d}^* & (1- b^2 u^2 ) {\bf I}_R  
\end{array}
\right]
.
\end{array}
\]
By 5,6 and 7 of Lemma 3, we have 
\[
\begin{array}{rcl}
{\bf ML} & = & 
\left[
\begin{array}{cc}
(1- b^2 u^2 ) {\bf I}_V & 0 \\ 
u(1- b^2 u^2 ) {\bf S} {\bf d}^* & -cu {\bf S} {\bf d}^* {\bf d} - bc u^2 {\bf S} {\bf d}^* {\bf d} {\bf S} +(1- b^2 u^2 ) {\bf I}_R 
\end{array}
\right]
\\ 
\  &   &                \\ 
\  & = & 
\left[
\begin{array}{cc}
(1- b^2 u^2 ) {\bf I}_V & 0 \\ 
u(1- b^2 u^2 ) {\bf S} {\bf d}^* & ( {\bf I}_R -u(c {\bf S} {\bf d}^* {\bf d} +b {\bf S} ))( {\bf I}_R +u b {\bf S} ) 
\end{array}
\right]
.
\end{array}
\]
Here, note that ${\bf S}^2 = {\bf I}_R$.

For $\mid t \mid , \mid u \mid $ sufficiently small, we have 
\[
\sigma ( \Delta (u)), \sigma ((1- b^2 t^2 ) {\bf I}_V ), 
\sigma ((1- b^2 t^2 ) {\bf I}_R ), 
\sigma (( {\bf I}_R -u(c {\bf S} {\bf d}^* {\bf d} +b {\bf S} ))( {\bf I}_R +u b {\bf S} ) )  
\]
\[
\in B_1 (1) = \{ z \in {\bf C} \mid \  \mid z-1 \mid <1 \} . 
\]
Similar to the proof of [9, Proposition 3.8], 
$\sigma ({\bf LM} )$ and $\sigma ({\bf ML})$ are contained in $B_1 (1)$. 
Thus, ${\bf L} $ and ${\bf M} $ are invertible, with bounded inverse, for 
$\mid t \mid , \mid u \mid $ sufficiently small. 

By 1 of Proposition 1, 1 of Proposition 2 and Corollary 1, we have 
\[
\begin{array}{rcl}
\det {}_{\Gamma } ( {\bf LM} ) & = & 
\det {}_{\Gamma } ((1- b^2 u^2 ) {\bf I}_V -cu {\bf d} {\bf S} {\bf d}^* -bc u^2 {\bf d} {\bf S}^2 {\bf d}^* ) 
\det {}_{\Gamma } ((1- b^2 u^2 ) {\bf I}_R) \\
\  &   &                \\ 
\  & = & 
(1- b^2 u^2 )^{{\rm Tr}_{\Gamma } ({\bf I}_R)} 
\det {}_{\Gamma } ((1-ab u^2 ) {\bf I}_V -cu {\bf d} {\bf S} {\bf d}^* ) 
\end{array}
\]
and 
\[
\begin{array}{rcl}
\det {}_{\Gamma } ( {\bf ML} ) & = & 
\det {}_{\Gamma } ((1- b^2 u^2 ) {\bf I}_V )
\det {}_{\Gamma } ( {\bf I}_R -u(c {\bf S} {\bf d}^* {\bf d} +b {\bf S} )) \det {}_{\Gamma } ( {\bf I}_R +u b {\bf S} ) \\
\  &   &                \\ 
\  & = & (1- b^2 u^2 )^{{\rm Tr}_{\Gamma } ({\bf I}_V)} 
\det {}_{\Gamma } ( {\bf I}_R -u(c {\bf S} {\bf d}^* {\bf d} +b {\bf S} )) \det {}_{\Gamma } ( {\bf I}_R +u b {\bf S} ) . 
\end{array}
\]

Let an orientation of $G$ be a choice of one oriented edge for each pair of 
edges in $R(G)$, 
which is called positively oriented. 
We denote by $E^+ G$ the set of positively oriented edges. 
Moreover, let $E^- G := \{ e^{-1} \mid e \in E^+ G \} $. 
An element of $E^- G$ is called a negatively oriented. 
Note that $R(G)=E^+ G \cup E^- G$. 

The operator ${\bf S}$ maps $\ell {}^2 (E^+ G)$ to $\ell {}^2 (E^- G)$. 
Then we obtain a representation $\rho $ of ${\cal B}( \ell {}^2 (R(G)))$ 
onto $Mat {}_2 {\cal B} ( \ell {}^2 (E^+ G))$, under 
\[
\rho ({\bf S} )=
\left[
\begin{array}{cc}
0 & {\bf I} \\
{\bf I} & 0 
\end{array}
\right]
, 
\rho ({\bf I}_R )=
\left[
\begin{array}{cc}
{\bf I} & 0 \\
0 & {\bf I} 
\end{array}
\right]
.
\]
By 1 and 3 of Proposition 2, 
\[
\begin{array}{rcl}
\det {}_{\Gamma } ({\bf I}_R +bu {\bf S}) 
& = & 
\det {}_{\Gamma }
\left[
\begin{array}{cc}
{\bf I} & -bu {\bf I} \\
0 & {\bf I} 
\end{array}
\right]
\  
\det {}_{\Gamma } 
\left[
\begin{array}{cc}
{\bf I} & bu {\bf I} \\
bu {\bf I} & {\bf I} 
\end{array}
\right]
\\ 
\  &   &                \\ 
\  & = & 
\det {}_{\Gamma } 
\left[
\begin{array}{cc}
(1- b^2 u^2 ){\bf I} & 0 \\
\ast & {\bf I} 
\end{array}
\right]
=(1- b^2 u^2 )^{\frac{1}{2} {\rm Tr}_{\Gamma} ( {\bf I}_R)} . 
\end{array}
\]

For $\mid t \mid , \mid u \mid $ sufficiently small, we have 
\[ 
{\bf ML}={\bf MLM}{\bf M}^{-1} ,
\]
and so, by 2 of Proposition 2, 
\[ 
\det {}_{\Gamma } ( {\bf LM} ) = \det {}_{\Gamma } ( {\bf ML} ) . 
\]

Therefore, it follows that 
\[
(1- b^2 u^2 )^{{\rm Tr}_{\Gamma } ({\bf I}_R)} 
\det {}_{\Gamma } ( (1-ab u^2 ) {\bf I}_V -cu {\bf d} {\bf S} {\bf d}^* )
\]
\[
= (1- b^2 u^2 )^{\frac{1}{2} {\rm Tr}_{\Gamma } ({\bf I}_R)
+ {\rm Tr}_{\Gamma } ({\bf I}_V)} 
\det {}_{\Gamma } ( {\bf I}_R -u {\bf S} (c {\bf d}^* {\bf d} +b {\bf I}_R )) , 
\]
and so 
\[
\begin{array}{rcl}
\  &  & \det {}_{\Gamma } ( {\bf I}_R -u {\bf S} {\bf C} )= \det {}_{\Gamma } ( {\bf I}_R -u {\bf S} (c {\bf d}^* {\bf d} +b {\bf I}_R )) \\
\  &   &                \\ 
\  & = & (1- b^2 u^2 )^{\frac{1}{2} {\rm Tr}_{\Gamma } ({\bf I}_R)
- {\rm Tr}_{\Gamma } ({\bf I}_V)} 
\det {}_{\Gamma } (  (1-ab u^2 ) {\bf I}_V -cu {\bf d} {\bf S} {\bf d}^* ) . 
\end{array}
\]

Hence the result follows by the definition of ${\rm TR}_{\Gamma } $.  
$\Box$

\end{document}